\documentclass{article}

\newtheorem{theorem}{Theorem}[section]
\newtheorem{corollary}[theorem]{Corollary}

\newtheorem{lemma}[theorem]{Lemma}
\newtheorem{proposition}[theorem]{Proposition}

\newtheorem{remark}{Remark}[section]

\evensidemargin\oddsidemargin
\begin{document}

\author{Vadim E. Levit and Eugen Mandrescu \\
Department of Computer Science\\
Holon Academic Institute of Technology\\
52 Golomb Str., P.O. Box 305\\
Holon 58102, ISRAEL\\
\{levitv, eugen\_m\}@barley.cteh.ac.il}
\title{The Intersection of All Maximum Stable Sets of a Tree and its Pendant
Vertices}
\date{}
\maketitle

\begin{abstract}
A \textit{stable} set in a graph $G$ is a set of mutually non-adjacent
vertices, $\alpha (G)$ is the size of a maximum stable set of $G$, and $%
core(G)$ is the intersection of all its maximum stable sets. In this paper
we demonstrate that in a tree $T$, of order $n\geq 2$, any stable set of
size $\geq n/2$ contains at least one pendant vertex. Hence, we deduce that
any maximum stable set in a tree contains at least one pendant vertex. Our
main finding is the theorem claiming that if $T$ does not own a perfect
matching, then at least two pendant vertices an even distance apart belong
to $core(T)$. While it is known that if $G$ is a connected bipartite graph
of order $n$ $\geq 2$, then $\left| core(G)\right| \neq 1$ (see Levit,
Mandrescu \cite{LeviMan2}), our new statement reveals an additional
structure of the intersection of all maximum stable sets of a tree. The
above assertions give refining of one result of Hammer, Hansen and Simeone 
\cite{HamHanSim}, stating that if a graph $G$ is of order less than $2\alpha
(G)$, then $core(G)$ is non-empty, and also of a result of Jamison \cite
{Jamison}, Gunter, Hartnel and Rall \cite{GunHarRall}, and Zito \cite{Zito},
saying that for a tree $T$ of order at least two, $\left| core(T)\right|
\neq 1$.
\end{abstract}

\section{Introduction}

Throughout this paper $G=(V,E)$ is a simple (i.e., a finite, undirected,
loopless and without multiple edges) graph with vertex set $V=V(G)$, edge
set $E=E(G)$, and its order is $\left| V\right| $. If $X\subset V$, then $%
G[X]$ is the subgraph of $G$ spanned by $X$. By $G-W$ we mean the subgraph $%
G[V-W]$, if $W\subset V(G)$. We also denote by $G-F$ the partial subgraph of 
$G$ obtained by deleting the edges of $F$, for $F\subset E(G)$, and we use $%
G-e$, if $F$ $=\{e\}$. Let $K_{n}$, $P_{n}$ denote the complete graph on $%
n\geq 1$ vertices and the chordless path on $n\geq 2$ vertices.

A set $A\subseteq V$ is \textit{stable} if no two vertices from $A$ are
adjacent. A stable set of maximum size will be referred as to a \textit{\
maximum stable set} of $G$, and the \textit{stability number }of $G$,
denoted by $\alpha (G)$, is the cardinality of a maximum stable set in $G$.
Let $\Omega (G)$ stand for the set $\{S:S$ \textit{is a maximum stable set of%
} $G\}$, and $core(G)=\cap \{S:S\in \Omega (G)\}$, (see \cite{LeviMan3}).

The \textit{neighborhood} of a vertex $v\in V$ is the set $N(v)=\{w:w\in V$
\ \textit{and} $vw\in E\}$, while the \textit{close neighborhood} of $v\in V$
is $N[v]=N(v)\cup \{v\}$. For $A\subset V$, we denote $N(A)=\{v\in
V-A:N(v)\cap A\neq \emptyset \}$, and $N[A]=N(A)\cup A$. In particular, if $%
\left| N(v)\right| =1$, then $v$ is a \textit{pendant vertex} of $G$. By $%
pend(G)$ we designate the set $\{v\in V(G):v$ \textit{is a pendant vertex in}
$G\}$.

By \textit{tree} we mean a connected acyclic graph of order greater than
one, and a \textit{forest} is a disjoint union of trees and isolated
vertices.

In this paper we show that any stable set $S$ of a tree $T$, of size $\left|
S\right| \geq \left| V(T)\right| /2$, contains at least one pendant vertex
of $T$. As a consequence, we infer that $S\cap pend(T)\neq \emptyset $ is
valid for any $S\in \Omega (T)$. Moreover, we prove that in a tree $T$ with $%
\alpha (T)>\left| V(T)\right| /2$, there exist at least two pendant vertices
belonging to every maximum stable set of $T$, such that the distance between
them is even.

We give also a new proof for a result of Hopkins and Staton, stating that if 
$\{A,B\}$ is the standard bipartition of the vertex set of a tree $T$, then $%
\Omega (T)=\{A\}$ or $\Omega (T)=\{B\}$ if and only if the distance between
any two pendant vertices of $T$ is even.

Our findings are also incorporated in the following contexts.

Firstly, the following theorem concerning maximum stable sets in general
graphs, due to Nemhauser and Trotter \cite{NemhTro}, shows that for a
special subgraph $H$ of a graph $G$, some maximum stable set of $H$ can be
enlarged to a maximum stable set of $G$. Namely, if $A\in \Omega (G[N[A]])$,
then there is $S\in \Omega (G)$, such that $A\subseteq S$. We show that, for
trees, some kind of an inverse theorem is also true. More precisely, we show
that any maximum stable set of a tree $T$ contains at least one of its
pendant vertices, i.e., for any $S\in \Omega (T)$ there exists some $A$,
such that $A\in \Omega (T[N[A]])$ and $A\subseteq S$.

Secondly, recall that Hammer, Hansen and Simeone have proved in \cite
{HamHanSim} that if a graph $G$ has $\alpha (G)>\left| V\left( G\right)
\right| /2$, then $\left| core(G)\right| \geq 1$. As a strengthening, Levit
and Mandrescu \cite{LeviMan2} showed that if $G$ is a connected bipartite
graph with $\left| V(G)\right| \geq 2$, then $\left| core(G)\right| \neq 1$.
Jamison \cite{Jamison}, Zito \cite{Zito}, and Gunther, Hartnel and Rall \cite
{GunHarRall} proved independently that $\left| core(T)\right| \neq 1$ is
true for any tree $T$. Now, for a tree $T$ with $\alpha (T)>\left|
V(T)\right| /2$, we demonstrate that $\left| pend(T)\cap core(T)\right| \geq
2$, which means that there exist at least two pendant vertices of $T$
belonging to all maximum stable sets of $T$.

Thirdly, it is well-known that any tree $T$ has at least two pendant
vertices (e.g., see Berge \cite{Berge}). Our results say that if $\alpha
(T)>\left| V(T)\right| /2$, then at least two pendant vertices of $T$ belong
to all maximum stable sets of $T$, and whenever $\alpha (T)=\left|
V(T)\right| /2$ then both parts of the standard bipartition of $T$ contain
at least one pendant vertex.\newpage

\section{Pendant vertices and maximum stable sets}

\begin{lemma}
\label{lem2}Any stable set consisting of only pendant vertices of a graph $G$
is contained in a maximum stable set of $G$.
\end{lemma}

\setlength {\parindent}{0.0cm}\textbf{Proof.} Let $A$ be a stable set of $G$
such that $A\subseteq pend(G)$, and $S\in \Omega (G)$. If $u\in A-S$, then $%
u $ is adjacent to some $w\in S-A$, otherwise $S\cup \left\{ u\right\} $ is
a stable set larger than $S$, which contradicts the maximality of $S$.
Hence, $S_{1}=S\cup \{u\}-\left\{ w\right\} \in \Omega (G)$, and $\left|
A\cap S\right| <\left| A\cap S_{1}\right| $. Therefore, using this exchange
procedure, after a finite number of steps, we have to obtain a maximum
stable set including $A$. \rule{2mm}{2mm}\setlength
{\parindent}{3.45ex}\newline

\begin{figure}[h]
\setlength{\unitlength}{1.0cm} 
\begin{picture}(5,1.2)\thicklines
  \put(4,0){\line(1,0){4}}
  \multiput(4,0)(2,0){3}{\circle*{0.25}}
  \multiput(5,0)(2,0){2}{\circle*{0.34}}
  \put(5,1){\circle*{0.25}}
  \put(4,0){\line(1,1){1}}
  \put(5,0){\line(0,1){1}}
  \put(6,1){\circle*{0.34}}
  \put(6,0){\line(0,1){1}}
  \put(5,1){\line(1,0){1}}
\end{picture}
\caption{$G$ has a maximum stable set containing no pendant vertex.}
\label{fig1}
\end{figure}
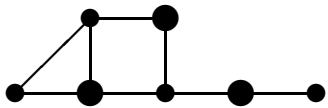

The converse of Lemma \ref{lem2} is not generally true. For instance, as it
is emphasized in Figure \ref{fig1}, the maximum stable set consisting of
only large vertices does not contain any pendant vertex of the graph $G$.

\begin{theorem}
\label{th4}If $S$ is a stable set of a tree $T$ and $\left| S\right| \geq
\left| V(T)\right| /2$, then $S\cap pend(T)$ is not empty. Moreover, there
exist $v\in S\cap pend(T)$ and $w\in S$ such that the distance between them
is $2$.
\end{theorem}

\setlength {\parindent}{0.0cm}\textbf{Proof.} Suppose, on the contrary, that 
$S\cap pend(T)=\emptyset $.\setlength
{\parindent}{3.45ex} Hence, any $s\in S$ has $\left| N(s)\right| \geq 2$.
Since $T$ is a tree and $S$ is a stable set of size $\left| S\right| \geq
\left| V(T)\right| /2$, it yields the following contradiction: 
\[
\left| V(T)\right| -1=\left| E(T)\right| \geq \left| (S,V(T)-S)\right| \geq
2\left| S\right| \geq \left| V(T)\right| . 
\]
Consequently, we infer that $S\cap pend(T)$ is not empty, for any stable set 
$S$ with $\left| S\right| \geq \left| V(T)\right| /2$.

We can assert now that there exists some $k\geq 1$, such that $S\cap
pend(T)=\{v_{i}:1\leq i\leq k\}$. Assume that for any $v_{i},w\in S$, the
distance between them is greater than two. If $N(v_{i})=\{u_{i}\},1\leq
i\leq k$, it follows that $S^{\prime }=(S-\{v_{i}:1\leq i\leq k\})\cup
\{u_{i}:1\leq i\leq k\}$ is a stable set in $T$ and $\left| S^{\prime
}\right| =\left| S\right| $. Clearly, $S^{\prime }\cap pend(T)=\emptyset $,
but this contradicts the fact that, according to the first part of the
theorem, $S^{\prime }\cap pend(T)$ must be non-empty, since $\left|
S^{\prime }\right| =\left| S\right| \geq \left| V(T)\right| /2$. Therefore,
there must exist $v\in S\cap pend(T),w\in S$ such that the distance between
them is $2$. \rule{2mm}{2mm}\newline

Now using the fact that $\alpha (G)\geq \left| V\left( G\right) \right| /2$
for any bipartite graph $G$, we obtain the following.

\begin{corollary}
\label{cor8}If $T$ is a tree, then $S\cap pend(T)\neq \emptyset $ for any $%
S\in \Omega (T)$.
\end{corollary}

\begin{figure}[h]
\setlength{\unitlength}{1.0cm} 
\begin{picture}(5,1.2)\thicklines
  \multiput(3,0)(1,0){4}{\circle*{0.34}}
  \multiput(3,1)(1,0){3}{\circle*{0.25}}
  \multiput(3,0)(1,0){3}{\line(0,1){1}}
  \multiput(3,0)(1,0){2}{\line(1,1){1}}
  \multiput(5,0)(1,0){2}{\line(-1,1){1}}
  \put(6,0){\line(-2,1){2}}

  \multiput(7,0)(2,0){2}{\circle*{0.34}}
  \put(8,0){\circle*{0.25}}
  \put(7,0){\line(1,0){2}}
  \multiput(7,1)(2,0){2}{\circle*{0.25}}
  \put(8,1){\circle*{0.34}}
  \put(8,1){\line(1,0){1}}
  \multiput(7,0)(1,0){3}{\line(0,1){1}}

\end{picture}
\caption{Bipartite graphs with maximum stable sets containing no pendant
vertices.}
\label{fig3}
\end{figure}
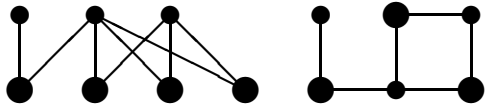

Corollary \ref{cor8} is not generally true for a connected graph $G$ with $%
pend(G)\neq \emptyset $ (see, for instance, the graph in Figure \ref{fig1}).
Notice also that it cannot be generalized to a bipartite graph $G$, both for 
$\alpha (G)=\left| V(G)\right| /2$ and $\alpha (G)>\left| V(G)\right| /2$
(see Figure \ref{fig3}). 

\begin{corollary}
\label{cor4}If $T$ is a tree with $\alpha (T)=\left| V(T)\right| /2$, and $%
\{A,B\}$ is its bipartition, then 
\[
A\cap pend(T)\neq \emptyset \ and\ B\cap pend(T)\neq \emptyset . 
\]
\end{corollary}

\setlength {\parindent}{0.0cm}\textbf{Proof.} In this case both $A$ and $B$
are maximum stable sets, because $A,B$ are stable and $\alpha (T)\geq \max
(\left| A\right| ,\left| B\right| )\geq \left| V\left( T\right) \right| /2,$%
and this implies $\left| A\right| =\left| B\right| =\alpha (T)$. Hence,
Corollary \ref{cor8} implies the result. \rule{2mm}{2mm}%
\setlength
{\parindent}{3.45ex}\newline

\begin{figure}[h]
\setlength{\unitlength}{1.0cm} 
\begin{picture}(5,1.2)\thicklines
  \multiput(5,1)(1,0){3}{\circle*{0.25}}
  \multiput(6,0)(1,0){2}{\circle*{0.25}}
  \put(6,0){\line(0,1){1}}
  \put(6,0){\line(-1,1){1}}
  \put(6,0){\line(1,1){1}}
  \put(7,0){\line(0,1){1}}
 \end{picture}
\caption{A tree with $A\cap pend(T)\neq \emptyset $, $B\cap pend(T)\neq
\emptyset $, and $\alpha (T)>\left| A\cup B\right| /2$.}
\label{fig4}
\end{figure}
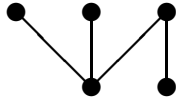

Figure \ref{fig4} shows that the converse of Corollary \ref{cor4} is not
generally true.

Since the distance between any two vertices belonging respectively to $A$
and $B$ is odd, we obtain the following form of Corollary \ref{cor4}.

\begin{corollary}
\label{cor5}If $T$ is a tree with $\alpha (T)=\left| V(T)\right| /2$, then $%
T $ contains at least two pendant vertices, such that the distance between
them is odd.
\end{corollary}

Recall from \cite{HopStat} that $G$ is a \textit{strong unique independent
graph} if $\left| \Omega (G)\right| =\left| \{S\}\right| =1$ and $V(G)-S$ is
also stable. For example, every cordless path of odd order belongs to this
class of graphs. Any strong unique independent graph $G$ is necessarily
bipartite, and its bipartition is $\{S,V(G)-S\}$. Using Theorem \ref{th4},
we are giving now an alternative proof of the following theorem
characterizing strong unique independent trees, which was first proved in 
\cite{HopStat}.

\begin{theorem}
\label{prop6}\cite{HopStat} If $\{A,B\}$ is the bipartition of the tree $T$,
then the following assertions are equivalent:

($\mathit{i}$) $T$ is a strong unique independent tree;

($\mathit{ii}$) $pend(T)\subseteq A$ or $pend(T)\subseteq B$;

($\mathit{iii}$) the distance between any two pendant vertices of $T$ is
even.
\end{theorem}

\setlength {\parindent}{0.0cm}\textbf{Proof.} ($\mathit{i}$) $\Rightarrow $ (%
$\mathit{ii}$) If $\Omega (T)=\{A\}$, then Lemma \ref{lem2} implies $%
pend(T)\subseteq A$.\setlength
{\parindent}{3.45ex}

The equivalence ($\mathit{ii}$) $\Leftrightarrow $ ($\mathit{iii}$) is clear.

($\mathit{ii}$) $\Rightarrow $ ($\mathit{i}$) Without loss of generality, we
may suppose that $\left| A\right| \geq \left| V(T)\right| /2$. Since $A$ is
also a stable set, Theorem \ref{th4} ensures that $A\cap pend(T)\neq
\emptyset $, and consequently, $pend(T)\subseteq A$. Let $S\in \Omega (T)$
and $S-A\neq \emptyset $. Then 
\[
\left| S\right| =\left| S\cap A\right| +\left| S\cap B\right| \geq \left|
S\cap A\right| +\left| N(S\cap B)\right| =\left| A\right| 
\]
and hence, $\left| S\cap B\right| \geq \left| N(S\cap B)\right| $. Since no
vertex in $S\cap A$ is adjacent to any vertex in $S\cap B$, it follows that
at least one tree, say $T^{\prime }$, of the forest $T[N[S\cap B]]$, has $%
\left| V(T^{\prime })\cap B\right| \geq \left| V(T^{\prime })\cap A\right| $%
. Consequently, by Theorem \ref{th4}, it proves that $T^{\prime }$ has at
least one pendant vertex, say $v$, in $V(T^{\prime })\cap B$. Since $%
N(v,T)=N(v,T^{\prime })$, we infer that some pendant vertex of $T$ must be
in $S\cap B$, in contradiction with $pend(T)\subseteq A$. Therefore, there
is no maximum stable set $S$ in $T$, such that $S-A\neq \emptyset $, and
since $A$ is a maximal stable set, it follows that, in fact, $\Omega
(T)=\{A\}$, i.e., $T$ is a strong unique independent tree. \rule{2mm}{2mm}%
\setlength
{\parindent}{3.45ex}

\begin{corollary}
If a tree $T$ has the bipartition $\{A,B\}$ and $S$ is a maximal stable set
such that $\left| S\right| >\min \{\left| A\right| ,\left| B\right| \}$,
then $S\cap pend(T)\neq \emptyset $.\label{cor1} Moreover, there exist $v\in
S\cap pend(T)$ and $w\in S$ such that the distance between them is $2$.
\end{corollary}

\setlength {\parindent}{0.0cm}\textbf{Proof.} Equivalently, we prove that if 
$S$ is a maximal stable set $S$ of $T$ satisfying $S\cap pend(T)=\emptyset $%
, then $\min \{\left| A\right| ,\left| B\right| \}\geq \left| S\right| $. If 
$S_{A}=S\cap A,S_{B}=S\cap B$, then $\{S_{A},B-S_{B}\}$ is the bipartition
of the tree $T_{1}=T[S_{A}\cup (B-S_{B})]$ and $\{A-S_{A},S_{B}\} $ is the
bipartition of the tree 
\[
T_{2}=T[(A-S_{A})\cup S_{B}]. 
\]
Since $(S_{A},S_{B})=\emptyset $, it follows that $S_{A}\cap
pend(T_{1})=\emptyset =S_{B}\cap pend(T_{2})$, and consequently, by Theorem 
\ref{prop6}, $T_{1}$ and $T_{2}$ are strong unique independent trees.
Therefore, both $\left| S_{B}\right| \leq \left| A-S_{A}\right| $ and $%
\left| S_{A}\right| \leq \left| B-S_{B}\right| $. Hence, we get that 
\[
\left| S\right| =\left| S_{A}\right| +\left| S_{B}\right| \leq \min \{\left|
S_{A}\right| +\left| A-S_{A}\right| ,\left| B-S_{B}\right| +\left|
S_{B}\right| \}=\min \{\left| A\right| ,\left| B\right| \}, 
\]
which completes the proof. \rule{2mm}{2mm}\setlength {\parindent}{3.45ex}%
\newline

\begin{figure}[h]
\setlength{\unitlength}{1.0cm} 
\begin{picture}(5,1.2)\thicklines

  \multiput(3,1)(1,0){5}{\circle*{0.25}}
  \put(4,0){\circle*{0.25}}
  \multiput(6,0)(1,0){3}{\circle*{0.25}}
  \put(4,0){\line(-1,1){1}}
  \put(4,0){\line(0,1){1}}
  \put(4,0){\line(1,1){1}}
  \put(4,0){\line(2,1){2}}
  \put(6,0){\line(0,1){1}}
  \put(6,0){\line(1,1){1}}
  \put(7,0){\line(0,1){1}}
  \put(8,0){\line(-1,1){1}}
 
\put(3.6,0){\makebox(0,0){$a$}}
\put(7.3,1){\makebox(0,0){$b$}}

  \multiput(10,0)(1,0){2}{\circle*{0.25}}
  \multiput(10,1)(1,0){3}{\circle*{0.25}}
  \put(10,0){\line(0,1){1}}
  \put(10,0){\line(1,1){1}}
  \put(11,0){\line(0,1){1}}
  \put(11,0){\line(1,1){1}}
  \put(9.7,0){\makebox(0,0){$a$}}
  \put(10.7,0){\makebox(0,0){$b$}}

\end{picture}
\caption{$\{a,b\}$ is a maximal stabe set containing no pendant vertices.}
\label{fig3333}
\end{figure}
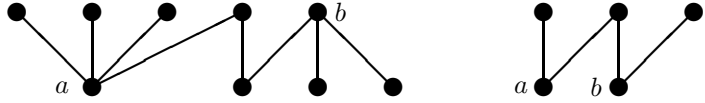

In other words, Corollary \ref{cor1} shows that for a maximal stable set $S$
of a tree with the bipartition $\{A,B\}$, it is enough to require that there
are no pendant vertices belonging to $S$ to ensure that $\min \{\left|
A\right| ,\left| B\right| \}\geq \left| S\right| $ (Figure \ref{fig3333}
shows examples of trees with $\min \{\left| A\right| ,\left| B\right|
\}=\left| S\right| $ and $\min \{\left| A\right| ,\left| B\right| \}>\left|
S\right| $). 

\begin{figure}[h]
\setlength{\unitlength}{1.0cm} 
\begin{picture}(5,1.2)\thicklines

  \multiput(5,1)(1,0){5}{\circle*{0.25}}
  \multiput(6,0)(2,0){2}{\circle*{0.25}}
  \put(6,0){\line(-1,1){1}}
  \put(6,0){\line(0,1){1}}
  \put(6,0){\line(1,1){1}}
  \put(8,0){\line(-1,1){1}}
  \put(8,0){\line(0,1){1}}
  \put(8,0){\line(1,1){1}}
   
\put(5.7,0){\makebox(0,0){$a$}}
\put(7.7,1){\makebox(0,0){$b$}}
\put(8.7,1){\makebox(0,0){$c$}}
 
\end{picture}
\caption{$\{a,b,c\}$ is a maximal stabe set of size $<n/2$ containing
pendant vertices.}
\label{fig3434}
\end{figure}
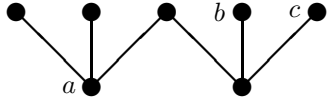

If $\left| A\right| \neq \left| B\right| $ then the claim of Corollary \ref
{cor1} is stronger than the corresponding direct consequence from Theorem 
\ref{th4}, because there is a tree $T$ containing a maximal stable set $S$,
such that $\min \{\left| A\right| ,\left| B\right| \}<\left| S\right|
<\left| V(T)\right| /2$ and $S\cap pend(T)\neq \emptyset $ (for an example,
see Figure \ref{fig3434}). 

\begin{figure}[h]
\setlength{\unitlength}{1.0cm} 
\begin{picture}(5,1.2)\thicklines

  \multiput(4,1)(1,0){5}{\circle*{0.25}}
  \multiput(5,0)(1,0){5}{\circle*{0.25}}
  \put(5,0){\line(-1,1){1}}
  \put(5,0){\line(0,1){1}}
  \put(5,0){\line(1,1){1}}
  \put(6,0){\line(0,1){1}}
  \put(6,0){\line(1,1){1}}
  \put(7,0){\line(0,1){1}}
  \put(8,0){\line(-1,1){1}}
  \put(8,0){\line(0,1){1}}
  \put(8,1){\line(1,-1){1}} 
\put(4.7,0){\makebox(0,0){$a$}}
\put(6.7,1){\makebox(0,0){$b$}}
\put(9.3,0){\makebox(0,0){$c$}}
 
\end{picture}
\caption{$\{a,b,c\}$ is a maximal stabe set of size $<\min \{\left| A\right|
,\left| B\right| \}$ containing pendant vertices.}
\label{fig4545}
\end{figure}
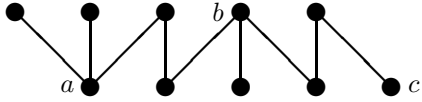

The converse of Corollary \ref{cor1} is not true, see, for instance, the
tree in Figure \ref{fig4545}. 

\section{Pendant vertices and intersection of all maximum stable sets}

We start this section with two results concerning $\alpha ^{+}$-stable
graphs, which we shall use in the sequel. Recall that a graph $G$ is $\alpha
^{+}$-\textit{stable} if $\alpha (G+e)=\alpha (G)$, for any edge $e\in E(%
\overline{G})$, where $\overline{G}$ is the complement of $G$, (see \cite
{GunHarRall}). The class of $\alpha ^{+}$-stable graphs was characterized by
Haynes et al. as follows:

\begin{theorem}
\label{th1}\cite{HayLawBrDu} A graph $G$ is $\alpha ^{+}$-stable if and only
if $\left| core(G)\right| \leq 1$.
\end{theorem}

In \cite{LeviMan} it was shown that:

\begin{theorem}
\label{th2}\cite{LeviMan} For a connected bipartite graph $G$ of order at
least two, the following assertions are equivalent:

($\mathit{i}$) $G$ is $\alpha ^{+}$-stable;

($\mathit{ii}$) $G$ has a perfect matching;

($\mathit{iii}$) $G$ possesses two maximum stable sets that partition its
vertex set;

($\mathit{iv}$) $\left| core(G)\right| =0.$
\end{theorem}

This statement generalizes the corresponding theorem of Gunther et al.,
proved for trees in \cite{GunHarRall}.

The following proposition and corollary have been proved for bipartite
graphs in the preprint \cite{LeviMan3}. Trying to give a self-consistent
presentation of this paper we prove them independently here.

\begin{proposition}
\label{prop4}If $T$ is a tree, then $\alpha (T)>\left| V\left( T\right)
\right| /2$ if and only if $\left| core(T)\right| \geq 2$.
\end{proposition}

\setlength {\parindent}{0.0cm}\textbf{Proof.} Suppose, on the contrary, that 
$\left| core(T)\right| \leq 1$. According to Theorem \ref{th1}, $T$ is an $%
\alpha ^{+}$-stable graph. Hence, Theorem \ref{th2} implies that $\alpha
(T)=\left| V\left( T\right) \right| /2$, in contradiction with the premise
on $\alpha (T)$.\setlength {\parindent}{3.45ex}

Conversely, if $\left| core(T)\right| \geq 2$, then Theorem \ref{th1}
ensures that $T$ is not $\alpha ^{+}$-stable. Since for trees $\alpha
(T)\geq \left| V\left( T\right) \right| /2$, Theorem \ref{th2} implies that $%
\alpha (T)>\left| V\left( T\right) \right| /2$. \rule{2mm}{2mm}

\begin{corollary}
\label{cor3}If $T$ is a tree, then $\alpha (T)=\left| V\left( T\right)
\right| /2$ if and only if $\left| core(T)\right| =0$.
\end{corollary}

Let $G_{i}=(V_{i},E_{i}),i=1,2$, be two graphs with $V_{1}\cap
V_{2}=\emptyset $, and $Q_{1},Q_{2}$ be cliques of the same size in $%
G_{1},G_{2}$, respectively. The \textit{clique bonding} of the graphs $%
G_{1},G_{2}$ is the graph $G=G_{1}*Q*G_{2}$ obtained by identifying $Q_{1}$
and $Q_{2}$ into a single clique $Q$, \cite{Berge}. In other words, $G$ is
defined by $V(G)=V_{1}\cup V_{2}-V(Q_{2})$ and

\[
E(G)=E_{1}\cup E_{2}-(\{xy:x,y\in V(Q_{2})\}\cup \{xy:x\in V(Q_{2}),y\in
V_{2}-V(Q_{2})\}). 
\]

If $V(Q)=\{v\}$, we shall denote the clique bonding of $G_{1}$ and $G_{2}$
by $G_{1}*v*G_{2}$.

\begin{lemma}
\label{lem3}If $T_{1},T_{2}$ are trees, $T=T_{1}*v*T_{2}$, and $v\in core(T)$%
, then 
\[
\alpha (T)=\alpha (T_{1})+\alpha (T_{2})-1. 
\]
\end{lemma}

\setlength {\parindent}{0.0cm}\textbf{Proof.} Let $S\in \Omega (T)$. Then $%
S\cap V(T_{i})$ is stable in $T_{i}$, and, therefore, it follows that $%
\left| S\cap V(T_{i})\right| \leq \alpha (T_{i})$, for each $i=1,2$. 
\setlength
{\parindent}{3.45ex}Hence, we get that 
\[
\alpha (T)=\left| S\right| =\left| S\cap V(T_{1})\right| +\left| S\cap
V(T_{2})-\{v\}\right| \leq \alpha (T_{1})+\alpha (T_{2})-1. 
\]

\textit{Case }$\mathit{1}$. There are $S_{i}\in \Omega (T_{i}),i=1,2$, such
that $v\notin S_{1}\cup S_{2}$. Hence, $S_{1}\cup S_{2}$ is stable in $T$
and $\left| S_{1}\cup S_{2}\right| =\alpha (T_{1})+\alpha (T_{2})>\alpha
(T_{1})+\alpha (T_{2})-1\geq \alpha (T)$, which brings a contradiction.

\textit{Case }$\mathit{2}$. There are $S_{i}\in \Omega (T_{i}),i=1,2$, such
that $v\in S_{1}\cap S_{2}$. Then $S_{1}\cup S_{2}$ is stable in $T$ and $%
\left| S_{1}\cup S_{2}\right| =\alpha (T_{1})+\alpha (T_{2})-1\geq \alpha
(T) $, and this implies that $\alpha (T)=\alpha (T_{1})+\alpha (T_{2})-1$.

\textit{Case }$\mathit{3}$. There are $S_{i}\in \Omega (T_{i})$ such that $%
v\in S_{1}$, and $v\notin S_{2}$. Hence, $(S_{1}-\{v\})\cup S_{2}$ is stable
in $T$ and $\left| (S_{1}-\{v\})\cup S_{2}\right| =\alpha (T_{1})+\alpha
(T_{2})-1\geq \alpha (T)$, and this assures that $\alpha (T)=\alpha
(T_{1})+\alpha (T_{2})-1$.

Consequently, we may conclude that $\alpha (T)=\alpha (T_{1})+\alpha
(T_{2})-1$. \rule{2mm}{2mm}

\begin{proposition}
\label{prop5}If $T_{1},T_{2}$ are trees, $T=T_{1}*v*T_{2}$, then $v\in
core(T)$ if and only if $v\in core(T_{i}),$ $i=1,2$.
\end{proposition}

\setlength {\parindent}{0.0cm}\textbf{Proof.} If $v\in core(T)$, then Lemma 
\ref{lem3} implies that $\alpha (T)=\alpha (T_{1})+\alpha (T_{2})-1$.
Suppose, on the contrary, that $v\notin core(T_{1})$. Let $S_{i}\in \Omega
(T_{i})$, $i=1,2$, be such that $v\notin S_{1}$. Hence $S=S_{1}\cup
S_{2}-\{v\}$ is stable in $T$ and $\left| S\right| =\alpha (T_{1})+\alpha
(T_{2})-1$. Consequently, $S\in \Omega (T)$ but $v\notin S$, in
contradiction with $v\in core(T)$.\setlength {\parindent}{3.45ex}

Conversely, let $v\in core(T_{i}),i=1,2$, $S_{i}\in \Omega (T_{i})$, $i=1,2$%
, and $S\in \Omega (T)$ be such that $v\notin S$. Then $S_{1}\cup S_{2}$ is
stable in $T$ and $\left| S_{1}\cup S_{2}\right| =\alpha (T_{1})+\alpha
(T_{2})-1$. Clearly, $S\cap V(T_{i})$ is stable in $T_{i},i=1,2$, and
because $v\in core(T_{i}),i=1,2$, we have that $\left| S\cap V(T_{i})\right|
\leq \alpha (T_{i})-1,i=1,2$. Hence, 
\[
\left| S\right| =\left| S\cap V(T_{1})\right| +\left| S\cap V(T_{2})\right|
\leq \alpha (T_{1})+\alpha (T_{2})-2<\left| S_{1}\cup S_{2}\right| , 
\]
and this contradicts the choice $S\in \Omega (T)$. \rule{2mm}{2mm}

\begin{lemma}
\label{lem4}If $T_{1},T_{2}$ are trees, $T=T_{1}*v*T_{2}$ and $v\in core(T)$%
, then 
\[
core(T)=core(T_{1})\cup core(T_{2}). 
\]
\end{lemma}

\setlength {\parindent}{0.0cm}\textbf{Proof.} According to Lemma \ref{lem3},
we have that $\alpha (T)=\alpha (T_{1})+\alpha (T_{2})-1$, and Proposition 
\ref{prop5} ensures that $v\in core(T_{i}),i=1,2$.%
\setlength
{\parindent}{3.45ex}

Let $w\in (core(T)-\{v\})\cap V(T_{1})$ and $S_{i}\in \Omega (T_{i}),i=1,2$.
Then $S_{1}\cup S_{2}\in \Omega (T)$, and, therefore, $w\in S_{1}$. Since $%
S_{1}$ is an arbitrary set from $\Omega (T_{1})$, we get that 
\[
w\in core(T_{1})\subset core(T_{1})\cup core(T_{2}),\ i.e.,\
core(T)\subseteq core(T_{1})\cup core(T_{2}). 
\]

Conversely, let $w\in core(T_{1})-\{v\}$, and suppose there is $S\in \Omega
(T)$, such that $w\notin S$. Let us denote $S_{i}=S\cap V(T_{i})$, for $%
i=1,2 $. Since $w\notin S_{1}$, it follows that $\left| S_{1}-\{v\}\right|
\leq \alpha (T_{1})-2$. Hence, we get a contradiction: 
\[
\left| S\right| =\left| S_{1}-\{v\}\right| +\left| S_{2}\right| \leq \alpha
(T_{1})-2+\alpha (T_{2})<\alpha (T_{1})+\alpha (T_{2})-1=\alpha (T)=\left|
S\right| . 
\]
Consequently, $core(T_{1})\cup core(T_{2})\subseteq core(T)$ is also valid,
and this completes the proof. \rule{2mm}{2mm}\newline

In the following statement we are strengthening Corollary \ref{cor8} and
Proposition \ref{prop4}.

\begin{theorem}
\label{th6}If $T$ is a tree with $\alpha (T)>\left| V(T)\right| /2$, then 
\[
\left| core(T)\cap pend(T)\right| \geq 2. 
\]
\end{theorem}

\setlength {\parindent}{0.0cm}\textbf{Proof.} According to Proposition \ref
{prop4}, $\left| core(T)\right| \geq 2$. Since $T$ is a tree, it follows
that $\left| pend(T)\right| \geq 2$. To prove the theorem we use induction
on $n=\left| V(T)\right| $. The result is true for $n=3$. Let $T=(V,E)$ be a
tree with $n=\left| V\right| >3$, and suppose that the assertion is valid
for any tree with fewer number of vertices.\ If $core(T)=pend(T)$, the
result is clear. If $core(T)\neq pend(T)$, let $v\in core(T)-pend(T)$ and $%
T_{1},T_{2}$ be two trees such that $T=T_{1}*v*T_{2}$. A bipartition of $%
N(v) $ gives rise to a corresponding division of $T$ into $T_{1}$ and $T_{2}$%
.\emph{\ }According to Proposition \ref{prop5}, $v\in core(T_{i}),i=1,2$.
Hence, Proposition \ref{prop4} implies that $\alpha (T_{i})>\left|
V(T_{i})\right| /2,i=1,2$. By the induction hypothesis, each $T_{i}$ has at
least two pendant vertices belonging to $core(T_{i})$. Lemma \ref{lem4}
ensures that $core(T)=core(T_{1})\cup core(T_{2})$, and, therefore, $T$
itself has at least two pendant vertices in $core(T)$. \rule{2mm}{2mm}%
\setlength
{\parindent}{3.45ex}

\begin{corollary}
Let $T$ be a tree with $\alpha (T)>\left| V\left( T\right) \right| /2$, and $%
k\geq 2$. If there exists a vertex $v\in core(T)$ of degree greater or equal
to $2k$, then $\left| core(T)\cap pend(T)\right| \geq 2k$.
\end{corollary}

\setlength {\parindent}{0.0cm}\textbf{Proof.} Let us partition $N(v)$ into $%
k $ subsets $N_{i}(v),1\leq i\leq k$, each one having at least two vertices.
Then we can write $T$ as 
\[
T=(...((T_{1}*v*T_{2})*v*T_{3})...)*v*T_{k}, 
\]

where $T_{i}$ is the subtree of $T$ containing $N_{i}(v)$ as the
neighborhood of $v$. Hence, by Lemma \ref{lem4}, it follows 
\[
pend(T)=\cup \{pend(T_{i}):1\leq i\leq k\},core(T)=\cup \{core(T_{i}):1\leq
i\leq k\}. 
\]

According to Proposition \ref{prop5}, $v\in core(T_{i}),1\leq i\leq k$, and,
therefore, Theorem \ref{th6} implies that: 
\[
\left| core(T)\cap pend(T)\right| =\left| core(T_{1})\cap pend(T_{1})\right|
+...+\left| core(T_{k})\cap pend(T_{k})\right| \geq 2k, 
\]

and this completes the proof. \rule{2mm}{2mm}\setlength {\parindent}{3.45ex}

\begin{remark}
For every natural number $k$ there exists a tree $T$ with a vertex $v$ of
degree $k$ such that $v\in core(T)$. For instance, such a tree $T=(V,E)$ can
be defined as follows: $V=\{v\}\cup \{x_{i}:1\leq i\leq 2k\}$ and $%
E=\{vx_{i}:1\leq i\leq k\}\cup \{x_{i}x_{i+k}:1\leq i\leq k\}$.
\end{remark}

Combining Theorems \ref{th2}, \ref{th6} and Proposition \ref{prop4}, we
obtain the following characterization of trees having no perfect matchings.

\begin{theorem}
\label{th7}If $T$ is a tree of order $n$, then the following assertions are
equivalent:

($\mathit{i}$) $\alpha (T)>n/2$;

($\mathit{ii}$) $T$ has no perfect matching;

($\mathit{iii}$) $\left| core(T)\cap pend(T)\right| \geq 2$;

($\mathit{iv}$) $\left| core(T)\right| \neq 0$.
\end{theorem}

\begin{theorem}
If $T$ is a tree with $\alpha (T)>\left| V(T)\right| /2$, then for at least
two distinct vertices from $core(T)\cap pend(T)$ the distance between them
is even. Moreover, if the set $core(T)\cap pend(T)$ contains exactly two
vertices, then the distance between them never equals $4$.
\end{theorem}

\setlength {\parindent}{0.0cm}\textbf{Proof.} Let $\{A,B\}$ be the
bipartition of $T$ into the color classes. Notice that the distance between
two vertices is even if and only if they belong to the same color class of $%
T $.\setlength
{\parindent}{3.45ex}

To prove the theorem we use induction on $n=\left| V(T)\right| $. If $n=3$,
then $T=P_{3}$ and the assertion is true. Let now $T$ be a tree with $n\geq
4 $ vertices. By Theorem \ref{th7}, $\alpha (T)>n/2$ yields $\left|
core(T)\cap pend(T)\right| \geq 2$.

\textit{Case }$\mathit{1}$\textit{.} If $\left| core(T)\cap pend(T)\right|
\geq 3$, then 
\[
\min (\left| A\cap core(T)\cap pend(T)\right| ,\left| B\cap core(T)\cap
pend(T)\right| )>1. 
\]
Hence at least two vertices of $core(T)\cap pend(T)$ belong to one color
class, i.e., the distance between them is even.

\begin{figure}[h]
\setlength{\unitlength}{1.0cm} 
\begin{picture}(5,1.2)\thicklines

  \multiput(3,1)(1,0){3}{\circle*{0.25}}
  \multiput(3,0)(1,0){6}{\circle*{0.25}}
  \put(3,1){\line(1,0){2}}
  \put(3,0){\line(1,0){5}}
  \put(5,1){\line(1,-1){1}}
     
\put(2.7,1){\makebox(0,0){$u$}}
\put(2.7,0){\makebox(0,0){$v$}}
\put(4,1.3){\makebox(0,0){$x$}}
\put(4,0.35){\makebox(0,0){$y$}}
\put(5,1.3){\makebox(0,0){$c$}}
\put(5,0.4){\makebox(0,0){$d$}}
 
\end{picture}
\caption{A tree $T$ with $\left| core(T)\cap pend(T)\right| =2$.}
\label{fig4444}
\end{figure}
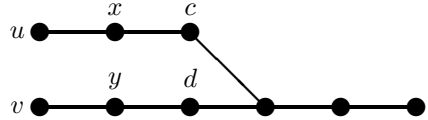

\textit{Case }$\mathit{2}$\textit{.} Let now $\left| core(T)\cap
pend(T)\right| =\left| \{u,v\}\right| =2$. Figure \ref{fig4444} shows that
such trees exist.

If $N(u)=N(v)$, then the distance between them is $2$, which is even.
Suppose now that $N(u)=\{x\}\neq N(v)=\{y\}$, and let $F=T[A\cup
B-\{u,v,x,y\}]$. Since $u$ and $v$ belong to all maximum stable sets of $T$,
we conclude that neither $x$ nor $y$ are contained in any maximum stable set
of $T$. Hence $\Omega (T)=\{S\cup \{u,v\}:S\in \Omega (F)\}$. Consequently, $%
core(T)=core(F)\cup \{u,v\}$ and $\alpha (F)=\alpha
(T)-2>n/2-2=(n-4)/2=\left| V(F)\right| /2$. Suppose that $F$ consists of $%
k\geq 1$ disjoint trees $\{T_{i}:1\leq i\leq k\}$. Since $\alpha (F)=\alpha
(T_{1})+...+\alpha (T_{k})>\left| V(F)\right| /2$, at least one tree, say $%
T_{j}$, has $\alpha (T_{j})>\left| V(T_{j})\right| /2$. By the induction
hypothesis, there exist $c,d\in core(T_{j})\cap pend(T_{j})$ such that the
distance between them in $T_{j}$ is even.

The pair of vertices $\{c,d\}\subset N(x)\cup N(y)$. Otherwise, if, for
instance, $c\notin N(x)\cup N(y)$, then $c\in core(T)\cap pend(T)$ and this
contradicts the fact that $core(T)\cap pend(T)=\{u,v\}$. If $\{c,d\}\subset
N(x)$ or $\{c,d\}\subset N(y)$, then $T_{j}$ is not a tree, since $\{cx,xd\}$
or $\{cx,xd\}$ builds a new path connecting $c$ and $d$ in addition to the
unique path between $c$ and $d$ in $T_{j}$ (together the two paths build a
cycle, which is forbidden in trees). Suppose that $\{c\}\subset N(x)$ and $%
\{d\}\subset N(y)$. If $xy\in $ $E(T)$ then again we see that $T_{j}$ can
not be a tree.

No edge from the set $\{uv,uy,ud,uc,vx,vc,vd\}$ exists since the vertices $u$
and $v$ are pendant in $T$. One can find an example of such a situation in
Figure \ref{fig4444}. The vertices $c$ and $d$ are not adjacent in $T$
because they are pendant in $T_{j}$. Thus, the distance between $u,v$ in $T$
is greater than the distance between $c$ and $d$ in $T_{j}$ by $4$, and
consequently, it is even.

If the distance between $u$ and $v$ is not equal to $2$, then $x\neq y$. Now
the same reasoning as above brings us to the conclusion that the condition $%
\left| core(T)\cap pend(T)\right| =2$ implies $\left| N(x)\cup N(y)\right|
\geq 4$, where the shortest path between $u$ and $v$ goes through the
vertices $x,c,y,d$, at least. Hence, the distance between $u$ and $v$ is
different from $4$. See Figure \ref{fig4444} for illustration of this claim. 
\rule{2mm}{2mm}

\section{Conclusions}

In this paper we have studied relationships between pendant vertices and
maximum stable sets of a tree. We have obtained a more precise version of
the well-known result of Berge, \cite{Berge}, stating that $\left|
pend(T)\right| \geq 2$ holds for any tree $T$ having at least two vertices.
Namely, we have proved that for such a tree $T$ either it has a perfect
matching and then both $A\cap pend(T)\neq \emptyset $ and $B\cap pend(T)\neq
\emptyset $, where $\{A,B\}$ is its bipartition, or it has not a perfect
matching and then at least two of its pendant vertices an even distance
apart belong to all maximum stable sets. As open problems, we suggest the
following: are there at least two pendant vertices of $T$ belonging to $\cap
\{S:S$ \textit{is a maximal stable set in }$T$\textit{\ of size} $k\}$, for $%
k=\left| V(T)\right| /2$, or $k=\min \{\left| A\right| ,\left| B\right| \}$?

\end{document}